\documentstyle[12pt]{article}
\begin{document}
 
\def\ba{\begin{eqnarray}}
\def\ea{\end{eqnarray}}
\def\lb{\label}
\def\ni{\noindent}
\def\nn{\nonumber}
\def\cl{\centerline}
\def\be{\begin{equation}}
\def\ee{\end{equation}}

\def\R{\hat{R}}
\def\F{\hat{F}}
\def\id{ I}
\def\M{{{{\cal M}}(\R,\F)}}

\def\s{\sigma}
\def\a{\alpha}
\def\b{\beta}
\def\g{\gamma}
\def\d{\delta}
\def\m{\mu}
\def\n{\nu}
\def\e{\epsilon}
\def\t{\tau}

\def\ok{\overline{k}}
\def\uk{\underline{k}}
\def\D{{{ D}}}

\def\rank{{\mbox{\rm rank}\,}}
\def\Aut{{\mbox{\rm Aut}\,}}
\def\End{{\mbox{\rm End}\,}}
\def\tr{{\mbox{\rm Tr}}\,}
\def\trq{{{\mbox{\rm Tr}}_{\scriptscriptstyle  \! \F}}}
\newcommand{\eod}{\hfill \rule{2.5mm}{2.5mm}}
\def\for{\mathrm{for~}} 

\title{ \vskip-3cm
\hfill{\normalsize\bf Preprint CPT-98/P.3700}\\
\vskip4cm
On quantum matrix algebras satisfying the 
Cayley-Hamilton-Newton identities}

\author{A. Isaev\\
{\small\em Bogoliubov Laboratory of Theoretical Physics, JINR, 141980  Dubna,}\\
{\small\em
 Moscow region, Russia} \medskip \\
O. Ogievetsky\footnote{On leave of absence from 
P. N. Lebedev Physical Institute, 
Theoretical Department, Leninsky pr. 53, 117924 Moscow, Russia}
~~~and~~~ P. Pyatov\footnote{On leave of absence from
Bogoliubov Laboratory of Theoretical Physics, JINR,
141980  Dubna, Moscow region, Russia}\\
{\small\em Center of Theoretical Physics, Luminy,
13288 Marseille, France}}
\date{}
\maketitle

\begin{abstract}
The Cayley-Hamilton-Newton identities which generalize both the characteristic
identity  and the Newton relations
have been recently obtained for the  algebras of the RTT-type.
We extend this result to a wider class of algebras $\M$
defined by a pair of compatible solutions of the Yang-Baxter equation.
This class includes the RTT-algebras as well as the Reflection equation algebras.
\end{abstract}
\newpage

In last years the two basic statements
of matrix algebra ---
the Cayley-Hamilton theorem and the Newton relations --- were 
generalized for quantum matrix algebras of the "RTT-" and the "Reflection equation"
(RE) types \cite{EOW,NT,PS,GPS,PS2,IOPS}. In   \cite{IOP} a new 
family of matrix identities called the Cayley-Hamilton-Newton (CHN) 
identities have been constructed.
The Cayley-Hamilton theorem and the Newton relations are  
particular cases and combinations of these identities. However the proof of the CHN
identities given in \cite{IOP} is adapted for the RTT-algebra case.
The factorization map from the RTT-algebra to the RE algebra produces,
in the quasitriangular case, the CHN identities for the RE algebra. 
In the present note we introduce a wider class of  algebras 
and extend for them the proof of the CHN identities given in \cite{IOP}.

The key observation for such a generalization is that there are
two R-matrices lying behind the construction of the CHN identities.
First of them which we denote  $\R$ 
is an R-matrix of the Hecke type. It is responsible, roughly speaking, for  
the commutation relations of quantum matrix entries. The second one which will be
referred as $\F$  is a closed R-matrix and it performs transition between
different matrix spaces. 
These two R-matrices are related by certain
compatibility conditions (see below, equations (\ref{compat})).

While the role of the first R-matrix is widely recognized, 
the importance of the second one is usually not noticed.
 In the case of the RTT-algebra the R-matrix $\F$ coincides with the
permutation matrix $P$ whereas for the RE algebra one has
$\F=\R$. Therefore $\F$ is in some sense trivial
for these standard examples of quantum matrix algebras.
Revealing  an independent role  of the R-matrix $\F$
allows to broaden the class of algebras under consideration 
and to give a universal proof of the CHN identities
for this whole class.
\vspace{2mm}

\ni {\bf 1. Notation.~}
Consider a pair of R-matrices $\R,  \F\in \Aut(V\otimes V)$ where 
$V$ is a finite-dimensional vector space. We call them {\em compatible}
if, besides the Yang-Baxter equations
\be
\lb{ybe}
\R_1\R_2\R_1=\R_2\R_1\R_2\,  , \qquad \F_1\F_2\F_1=\F_2\F_1\F_2\, ,
\ee
they satisfy the  conditions
\be
\lb{compat}
\R_1\F_2\F_1=\F_2\F_1\R_2\,  , \qquad \F_1\F_2\R_1=\R_2\F_1\F_2\, .
\ee
We use here  the matrix conventions of \cite{FRT}.
In particular, $\R_k$ and $\F_k$ denote the R-matrices acting 
in $V_k\otimes V_{k+1}$ --- the $k$-th and the $(k+1)$-st copies of the space $V$. 

In the sequel, we assume that $\R$ and $\F$ are compatible.  
Further on, we assume that $\R$
is an {\em even Hecke} R-matrix of a {\em height} $n$ and
$\F$ is a {\em closed} R-matrix.  Below we remind briefly these notions  
(for more details on the notation 
see \cite{G,GPS}). 

{\bf Conditions on the matrix $\R$.~} An R-matrix $\R$ 
satisfying the condition 
\be
\lb{hecke}
\R^2=\id +(q-q^{-1})\R \ 
\ee
is called the Hecke  R-matrix.
Here $\id$ is the identity operator and  $q\neq 0$ is a number.

Given a Hecke R-matrix, one constructs two sequences of projectors, 
$A^{(k)}$ and  $S^{(k)}\in \End(V^{\otimes k})$,
called {\em $q$-antisymmetrizers} and {\em $q$-symmetrizers}, correspondingly.
They are defined inductively, 
\ba
\lb{antis}
A^{(1)}:=\id\  ,&\quad&
A^{(k)}:={1\over k_q}\,
A^{(k{-}1)}\left(q^{k-1}-(k{-}1)_q\R_{k{-}1}\right)A^{(k{-}1)}\ ,
\\
\lb{simm}
S^{(1)}:=\id\  ,&\quad&
S^{(k)}:= {1\over k_q}\,
S^{(k{-}1)}\left(q^{1-k}+(k{-}1)_q\R_{k{-}1}\right)S^{(k{-}1)}\ ,
\ea
where it is additionally supposed  that
$k_q := (q^k - q^{-k})/(q-q^{-1}) \neq 0,~ \forall k=1,2,\dots~$. 

The Hecke R-matrix $\R$ is called even if its sequence of $q$-antisymmetrizers
vanishes at the $(n+1)$-st step and $~\rank A^{(n)} = 1$.
The number $n$ is called then the height of $\R$.
 
{\bf Conditions on the matrix $\F$.~} 
An R-matrix $\F=\F^{ab}_{cd}$ is called the closed R-matrix provided it is
invertible in indices $(a,c)$ and nonsingular (i.e., invertible in indices $(a,b)$). 
The first condition means that 
there exists a matrix $\Psi^{ab}_{cd}$
satisfying  $\Psi^{af}_{cg}\F^{gb}_{fd}=\d^a_d \d^b_c$ 
(summation over repeated indices is always assumed).
Denote $(\D)^a_b=\Psi^{ac}_{bc}$ .
Using the matrix $\D$, one introduces the notion of the {\em quantum trace}
for an arbitrary (not necessarily with commuting entries) matrix $X$,
\be
\trq X:=\tr\D X\ .
\lb{tr}
\ee

The following properties of the matrix $\D$ will be important for us
\ba
\lb{tr-f}
&&\trq_{(2)} \F_1 = \id_1\ ,  
\\[1mm]
\lb{fdd}
&&\F_1 \D_1 \D_2 = \D_1 \D_2 \F_1\ ,  
\\[1mm]
\lb{inv}
&&\trq_{(2)} \F_1^{\pm 1} \, X_1 \, {\F_1}^{\mp 1} = \id_1 \,
\trq X\ .  
\ea
Here and below we use notation
$\tr_{(i_1\dots i_k)}$
to denote the operation of taking traces in spaces 
with the numbers $i_1\dots i_k$.

{\bf Properties of compatible R-matrices.~}
Due to the compatibility conditions (\ref{compat}),
a matrix 
\be
\lb{twist}
\R^{\F} := \F \R \F^{-1}\ ,
\ee
satisfies the Yang-Baxter equation and is again compatible with $\F$.
This transformation was called {\em twisting} of R-matrices \cite{D}
and, in the case of compatible $\R$ and $\F$,
it has been considered in \cite{R}.

Since $\R^{\F}$ and $\F$ are compatible, one can consider the square of the twist, 
$\R^{\F\F}:=(\R^{\F})^{\F}$. One has the following relation
\be
\lb{twist2}
\R^{\F\F}_1 \D_1 \D_2 = \D_1 \D_2 \R_1\ .
\ee
We give the proof of this relation since 
we could not find it in literature.

Let $Y_{12}$ denote an arbitrary element of $\End(V\otimes V)$.  
Consider the following chain of transformations 
\ba
\nn
&&
\phantom{a}\!\!\!\!\!\!
\trq_{(1,2)} ( \R_1 \F_1^{2} Y_{12} ) =
\trq_{(1,2,3)} ( \R_1 \F_1 \F_2 \F_1 Y_{12} ) =        
\trq_{(1,2,3,4)} ( \R_1 \F_2 \F_1 \F_3 \F_2 Y_{12} ) 
\\[1mm]
&& =
\nn  
\trq_{(1,2,3,4)} ( \F_2 \F_1 \F_3 \F_2 Y_{12} \R_3) = 
\trq_{(1,2,3,4)} ( Y_{12} \R_3 \F_2 \F_3 \F_1 \F_2 ) 
\\[1mm]
&&=
\lb{transform}
\trq_{(1,2,3,4)} ( Y_{12} \F_2 \F_3 \F_1 \F_2 \R_1) =
\trq_{(1,2)} ( Y_{12} \F_1^{2} \R_1 ) \; .
\ea
Here we have subsequently used the
equations (\ref{tr-f}) and (\ref{ybe}) in the first line, (\ref{compat})
and the cyclic property of the trace together with (\ref{fdd}) in the second line,
and again (\ref{compat}), (\ref{tr-f}) and (\ref{ybe}) in the last line of the calculation.
Substituting the definition of the quantum trace, 
the result of (\ref{transform}) can be presented in a form
\[
\tr_{(1,2)}(Y_{12} \D_1 \D_2 \R_1 \F_1^2) = \tr_{(1,2)}(Y_{12} \F_1^2 \R_1 \D_1 \D_2)
\] 
which reduces to (\ref{twist2}) if one takes into account the arbitrariness of
$Y_{12}$ and applies once again the equation (\ref{fdd}).
\vspace{2mm}

\ni {\bf 2.  Algebra $\M$.~}
Consider a matrix $M$.  Usually one associates with $M$ a series
of its copies $M_k$  acting on the corresponding 
vector space $V_k$, $k=1,2,\dots$. 
We need the following generalization of this notion.

With a matrix $M$, we associate a series
of matrices $M_{\ok}$ defined inductively as
\be
\lb{Mk}
M_{\overline{1}} := M_1\ , \qquad M_{\overline{k+1}}:= \F_k M_{\ok} {\F_k}^{-1}\ .
\ee 
For $\F=P$ the new notation coincides with the old one:
$M_{\ok}\equiv M_k$. In general, the operator $M_{\ok}$ acts nontrivially on
the space $V_1\otimes\dots\otimes V_k$, not necessarily on $V_k$ alone.

Now we define the main object of this note, the algebra $\M$.
It is a unital associative algebra, generated by the components 
of  a matrix $M$ subject to a relation
\be
\lb{QMA}
\R_1 M_{\overline{1}} M_{\overline{2}} =
M_{\overline{1}} M_{\overline{2}} {\R^{\F\F}}_1\ ,
\ee
or 
$\R_1 M_1 \F_1 M_1 = M_1 \F_1 M_1 \R^F_1 $, 
in old notation.
Specializing to $\F=P$ or $\F=\R$ one reproduces the
RTT- or RE algebras, respectively. 
The algebras $\M$ form a subclass  
of more general algebras discussed in \cite{FM,KSk}.

In the Lemma below we collect several useful results.

\ni
{\bf Lemma.~}{\em
{\bf a)~} For a matrix $M$ with arbitrary entries, the following relations hold
\ba
\lb{l1}
&&\F_i M_{\ok} = M_{\ok} \F_i\ , \qquad \for  k\neq i, i+1 ,  
\\[1mm]
\lb{l1a}
&&\R_i M_{\ok}= M_{\ok} \R_i\ , \qquad \for k\neq i, i+1 ,
\\[1mm]
\lb{l2}
&&\F_{i\rightarrow k} M_{\overline{i}} M_{\overline{i+1}}\dots M_{\ok} =
M_{\overline{i+1}} M_{\overline{i+2}}\dots M_{\overline{k+1}}\F_{i\rightarrow k} ,
\quad \for \, i\leq k .
\ea
Here $\F_{i\rightarrow k}:= \F_i\F_{i+1}\dots\F_k$.
\vspace{1mm}

{\bf b)~}
Let $Y^{(k)}\equiv Y^{(k)}(\R_1,\dots ,\R_{k-1})$ be any polynomial
in $\R_1,\dots ,\R_{k-1}$, and let $Y^{(i,k)}:= Y^{(k)}(\R_i,\dots ,\R_{i+k-2})$.
Denote  $$\a(Y^{(k)}):=\trq_{(1,\dots ,k)}(Y^{(k)}M_{\overline{1}}\dots
M_{\ok})\ .$$

For a matrix $M$ with  arbitrary entries one has
\be
\lb{l4}
\trq_{(i,\dots ,i+k-1)}(Y^{(i,k)}M_{\overline{i}}\dots
M_{\overline{i+k-1}}) = \id_{1,\dots ,i-1}\, \a(Y^{(k)})\  ,
\ee
where $\id_{1,\dots ,i-1}$ is the identity in the spaces $1,\dots ,i-1$.

{\bf c)~} If, in addition, $M$ is the matrix of generators of $\M$,
one has 
\be
\lb{l3}
\R_k M_{\ok} M_{\overline{k+1}} = M_{\ok}
M_{\overline{k+1}} {\R^{\F\F}}_k\ .
\ee
}

\ni
{\bf Proof.~}
{\bf a)~} The relations (\ref{l1}) and (\ref{l1a})
are trivial for $i>k$. For $i<k-1$, the relations (\ref{l1}) and (\ref{l1a})
follow immediately from our definition (\ref{Mk}) of
$M_{\overline{i}}$ and the conditions (\ref{ybe}) and (\ref{compat}).

The relation (\ref{l2}) can be proved by induction.
For $k=i$ the formula (\ref{l2}) is just the definition of
$M_{\overline{k+1}}$. Suppose that (\ref{l2}) is valid for
some $k=j-1\geq i$. Then for $k=j$ we have
\ba
\nn
&&\F_{i\rightarrow j}M_{\overline{i}} 
M_{\overline{i+1}}\dots M_{\overline{j}} =
(\F_{i\rightarrow j-1}M_{\overline{i}}\dots 
M_{\overline{j-1}})(\F_j M_{\overline{j}}) 
\\[1mm]
&& =
M_{\overline{i+1}}\dots M_{\overline{j}}
(\F_{i\rightarrow j-1} M_{\overline{j+1}})\F_j =
M_{\overline{i+1}}\dots M_{\overline{j}} 
M_{\overline{j+1}}\F_{i\rightarrow j}\ ,
\nn
\ea
which completes the induction.
Here  we applied several times the  relations
(\ref{l1}), used the induction assumption and the definition of
$\F_{i\rightarrow j+1}$ and $M_{\overline{j+1}}$.

{\bf b)~} It suffices to check (\ref{l4}) for the case $i=2$. 
The calculation proceeds as follows
\ba
\nn
&&
\phantom{a}\!\!\!\!\!\!\!\!
\trq_{(2,\dots ,k+1)}(Y^{(2,k)}M_{\overline{2}}\dots M_{\overline{k+1}}) =
\trq_{(2,\dots ,k+1)}(Y^{(2,k)} \F_{1\rightarrow k}
M_{\overline{1}}\dots M_{\overline{k}}{\F_{1\rightarrow k}}^{-1}) 
\\[1mm]
\nn
&&
\phantom{a}\!\!\!\!\!\!\!\!
=
\trq_{(2,\dots ,k+1)} (\F_{1\rightarrow k}Y^{(k)}
M_{\overline{1}}\dots M_{\overline{k}}{\F_{1\rightarrow k}}^{-1}) 
\\[1mm]
\nn
&&\phantom{a}\!\!\!\!\!\!\!\!
=
\trq_{(2,\dots ,k)} (\F_{1\rightarrow k-1}\left[
\trq_{(k+1)}(\F_k Y^{(k)} M_{\overline{1}}\dots M_{\overline{k}}
{\F_k}^{-1})
\right]{\F_{1\rightarrow k-1}}^{-1}) 
\\[1mm]
\nn
&&
\phantom{a}\!\!\!\!\!\!\!\!
=
\trq_{(2,\dots ,k)} (\F_{1\rightarrow k-1}\left[
\id_k
\trq_{(k)}(Y^{(k)} M_{\overline{1}}\dots M_{\overline{k}})
\right]{\F_{1\rightarrow k-1}}^{-1}) =
\dots = \id_1\,  \a(Y^{(k)})  .
\ea
Here we used the equations 
(\ref{l2}), (\ref{compat}) and (\ref{inv}).
One should not be confused with the appearance of two 
$\trq_{(k)}$ in the left part of the last line of calculation.
The inner of these quantum traces acts on  arguments
in parentheses while the outer one respects only the
identity operator $\id_k$ among the terms enclosed by the square
brackets. Therefore the outer quantum trace $\trq_{(k)}$
can be calculated in the next step and transformed into an inner
$\trq_{(k-1)}$. The procedure repeats until all the outer
quantum traces transform into inner ones.

{\bf c)~} Induction in $k$. The relation (\ref{l3}) with $k=1$ is just the definition of $\M$. 
Assume that (\ref{l3}) is true
for some $k = i-1\geq 1$ and consider the case $k=i$,
\ba
\nn
&&\R_i M_{\overline{i}} M_{\overline{i+1}} =
\R_i \F_{i-1} M_{\overline{i-1}} M_{\overline{i+1}}{\F_{i-1}}^{-1} =
\R_i \F_{i-1} \F_i M_{\overline{i-1}} M_{\overline{i}}({\F_{i-1}\F_i})^{-1} 
\\[1mm]
\lb{proof}
&&=
\F_{i-1}\F_i M_{\overline{i-1}} M_{\overline{i}}({\F_{i-1}\F_i})^{-1}
{\R^{\F\F}}_i =
M_{\overline{i}} M_{\overline{i+1}} {\R^{\F\F}}_i\ .
\ea
Here we applied, first, the definition of $M_{\overline{i}}$, $M_{\overline{i+1}}$
and the  relations
(\ref{l1}). Next,  we used (\ref{compat}) and the induction assumption
and, then, performed the transformations of the first line 
of (\ref{proof}) in the inverse order.
\eod
\vspace{2mm}

\ni {\bf 3.  Characteristic subalgebra.~}
Let us consider three sequences of elements 
of the algebra $\M$:
\ba
\lb{s}
&&s_k(M) :=  
\trq_{(1\dots k)}(\R_{1\rightarrow k-1}M_{\overline{1}} M_{\overline{2}}
\dots M_{\ok})\ ,
\\[1mm]
\lb{sig}
&&\s_k(M) := 
\trq_{(1\dots k)}(A^{(k)}M_{\overline{1}} M_{\overline{2}}
\dots M_{\ok})\ ,
\\[1mm]
\lb{t}
&&\t_k(M) := 
\trq_{(1\dots k)}(S^{(k)}M_{\overline{1}} M_{\overline{2}}
\dots M_{\ok})\ ,
\quad k=1,2,\dots\; .
\ea
Also we put $s_0(M) =\s_0(M) = \t_0(M) = 1$.

These elements are interpreted as  symmetric polynomials
on the spectrum of the matrix $M$ (see \cite{IOPS,IOP}).
Namely, $s_k(M)$ are the {\em power sums}, $\s_k(M)$ are the
{\em elementary symmetric functions} and $\t_k(M)$ are the
{\em complete symmetric functions}. 

It follows from the Newton and Wronski relations (see below) that, given any pair of
the sets $\{s_k(M)\}$, $\{\s_k(M)\}$ or $\{\t_k(M)\}$, one can express
the elements of the first one of them as polynomials in the elements of the second one. 
Therefore all these sets generate the same subalgebra in $\M$ which
we call the {\em characteristic subalgebra} of $\M$.
\vspace{2mm}

\ni {\bf Proposition.~} 
{\em The characteristic subalgebra of $\M$ is  abelian.}
\vspace{2mm}

\ni{\bf Proof.~}
The commutativity of the characteristic subalgebra 
in the particular case of the RTT-algebra
was observed by J.M.Maillet \cite{M} who has 
checked the commutativity of power sums. 
We extend Maillet's method to treat the general case. 
The proof is based on the relation (\ref{l4}) which is 
trivial for the RTT-algebra  case but  
crucial for the general algebra $\M$. 

Consider a pair  $\a(Y^{(k)})$ and $\b(Z^{(i)})$ of elements of 
the characteristic subalgebra.
Using relations 
(\ref{l4}) one can present
the product of $\a$ and $\b$ in a form
\be
\lb{p1}
\a({Y^{(k)}})\,\b(Z^{(i)}) =
\trq_{(1,\dots ,k+i)}(Y^{(k)}Z^{(k+1,i)}M_{\overline{1}}
M_{\overline{2}}\dots M_{\overline{k+i}}) .
\ee

Further, consider an operator
$
U_{\scriptscriptstyle \R}:=
\R_{i\rightarrow i+k-1}\dots \R_{2\rightarrow k+1}\R_{1\rightarrow k}
$.
By virtue of the Yang-Baxter equation, one has
\be
\lb{shift}
Y^{(k)} = U_{\scriptscriptstyle \R}^{-1}Y^{(i+1,k)}U_{\scriptscriptstyle \R} \ , \qquad
Z^{(k+1,i)} = U_{\scriptscriptstyle \R}^{-1}Z^{(i)}U_{\scriptscriptstyle \R}\ .
\ee

Substituting (\ref{shift}) into (\ref{p1}), one continues the
transformation
\ba
\nn
&&\a({Y^{(k)}})\,\b(Z^{(i)}) =
\trq_{(1,\dots ,k+i)}(U_{\scriptscriptstyle \R}^{-1}Z^{(i)}Y^{(i+1,k)} U_{\scriptscriptstyle \R}
M_{\overline{1}}\dots M_{\overline{k+i}}) 
\\[1mm]
\nn
&&= 
\trq_{(1,\dots ,k+i)}(U_{\scriptscriptstyle \R}^{-1}
Z^{(i)}Y^{(i+1,k)}M_{\overline{1}}
\dots M_{\overline{i+k}}\ U_{\scriptscriptstyle \R^{\F\F}}) =
\b({Z^{(i)}})\,\a(Y^{(k)})\, .
\ea
Here the  relations (\ref{l1a}), (\ref{l3}), the cyclic property of the  trace
and the relation (\ref{twist2}) have been applied subsequently.
\eod
\vspace{2mm}

\ni {\bf 4.  Cayley-Hamilton-Newton identities.~}
Finally, we need a proper generalization of the notion of  a matrix power
for the case of  $\M$.
Taking off the first quantum trace in the definitions of symmetric
polynomials (\ref{s}), (\ref{sig}), (\ref{t}) one gets the following 
matrix expressions
\ba
\lb{power}
&&M^{\overline{k}}:= \trq_{(2,\dots ,k)}
(\R_{1\rightarrow k-1}M_{\overline{1}}M_{\overline{2}}\dots M_{\overline{k}})\ ,
\\[1mm]
\lb{skew-power}
&&M^{\wedge k}:= \trq_{(2,\dots ,k)}
(A^{(k)} M_{\overline{1}}M_{\overline{2}}\dots M_{\overline{k}})\ ,
\\[1mm]
\lb{simm-power}
&&M^{{\scriptscriptstyle\cal S} k}:= \trq_{(2,\dots ,k)}
(S^{(k)} M_{\overline{1}}M_{\overline{2}}\dots M_{\overline{k}})\ .
\ea
We call the matrix $M^{\overline{k}}$ the {\em $k$-th power
of the  matrix $M$}. Certainly, this definition coincides with
the usual one in the classical situation, $\R=\F=P$. 
More generally,  $M^{\overline{k}}\equiv M^k$ in the case $\R=\F$,
i.e., for the RE algebra.

The  matrices $M^{\wedge k}$ and $M^{{\scriptscriptstyle\cal S} k}$
will be relevant for the formulation of the Cayley-Hamilton-Newton identities.
It is natural to call them the {\em $k$-wedge} and the {\em $k$-symmetric} powers
of the  matrix $M$, respectively. 

With these definitions we can formulate our main result
\vspace{1mm}

\ni
{\bf Cayley-Hamilton-Newton theorem.~} 
{\em
Let  $M$ be the  matrix generating the algebra $\M$.
Then, the following matrix identities hold  
}
\ba
\lb{chn}
&&(-1)^{k-1}k_q M^{\wedge k} =
  \sum_{i=0}^{k-1}(-q)^i M^{\overline{k-i}}\ \s_i(M)\ ,
\\[1mm]
\lb{chn2}
&&k_q M^{{\scriptscriptstyle\cal S}k} =
  \sum_{i=0}^{k-1}q^{-i} M^{\overline{k-i}}\ \t_i(M)\ .
\ea
 
\ni{\bf Proof.~}
Consider the reflection 
$q\rightarrow -q^{-1}$, which is a symmetry
transformation of a parameter of the Hecke R-matrix $\R$.
It results in the substitutions 
$k_q\leftrightarrow (-1)^{k-1}k_q$, $A^{(k)}\leftrightarrow S^{(k)}$ 
and, hence, the two equations (\ref{chn}) and (\ref{chn2})
map into each other.
So, it suffices to prove  only one of the two series of 
equations (\ref{chn}) and (\ref{chn2}), say, the first one.

For the case of the RTT-algebra, these identities were proved in \cite{IOP}.
With the notation which we introduced in the present note,
the proof of these identities given in \cite{IOP}  
can be applied  practically without changes for the algebra 
$\M$.
The only additional remark should be given for the very first step of the proof.
It concerns the presentation of the typical term 
$M^{\overline{k-i}}\ \s_i(M)$ 
from the right hand side of the CHN identities in a form
\[
M^{\overline{k-i}} \s_i(M)
= \trq_{(2,\dots ,k)}
(\R_{1\rightarrow k-i-1}A^{(k-i+1,i)}
M_{\overline{1}}\dots M_{\overline{k}})\ .
\]
This equality being tautological in the RTT-algebra 
follows by an application of (\ref{l4}) and (\ref{l1a}) in the general case.

For the rest of the proof we refer the reader to the paper \cite{IOP}. \eod

In conclusion we present several corollaries of the Cayley-Hamilton-New\-ton
theorem. Their proofs given in \cite{IOP} for the case of the 
RTT-algebra remain valid for the general  algebra $\M$ as well. 

\ni
{\bf Newton relations.~}
\ba
\nn 
&&(-1)^{k-1}k_q \s_k(M) =
  \sum_{i=0}^{k-1}(-q)^i s_{k-i}(M) \s_i(M) ,
\\[1mm]
\nn
&&k_q \t_k(M) =
  \sum_{i=0}^{k-1}q^{-i} s_{k-i}(M) \t_i(M) .
\ea
\ni
{\bf Wronski relations.~}  
\ba
\nn
0 = \sum_{i=0}^{k}(-1)^{i} \t_{k-i}(M) \s_i(M) .
\ea
\ni
{\bf Cayley-Hamilton theorem.~} 
\ba
\nn
0 = \sum_{i=0}^n (-q)^i M^{\overline{n-i}} \s_i(M) ,
\quad
\mbox{where~~}M^{\overline{0}}:= q^{-n} n_q \tr_{(2,\dots ,n)}(A^{(n)}) \D^{-1} .
\ea
\ni
{\bf Inverse Cayley-Hamilton-Newton identities.~} 
\ba
\nn
M^{\overline{k}} = \sum_{i=1}^k (-1)^{i+1} q^{k-i} i_q M^{\wedge i} \t_{k-i}(M) =
\sum_{i=1}^k (-1)^{k-i} q^{i-k} i_q M^{{\scriptscriptstyle\cal S}i} \s_{k-i}(M) .
\ea
 
 \ni
 {\bf Acknowledgements}
 
 \ni
We are indebted to D.I. Gurevich and P.A. Saponov
for valuable discussions. This work was supported in parts
by  the CNRS grant PICS No.608 and the RFBR grant No. 98-01-2033.
The work of AI and PP  was also supported
by the RFBR grant No. 97-01-01041 and by the INTAS
grant 93-127-ext.  AI thanks CERN TH Division and OO thanks the Max-Planck-Institut
f\"ur Physik in M\"unchen for the hospitality.


\begin{thebibliography}{99}

\bibitem{EOW} Ewen H, Ogievetsky O and Wess J 1991 {\it Lett. Math.
Phys.} {\bf 22}  297--305

\bibitem{NT} Nazarov M and Tarasov V 1994 {\it  Publications RIMS} 
{\bf 30}  459--78

\bibitem{PS} Pyatov P N and Saponov P A 1995 {\it 
J. Phys. A: Math. Gen.} {\bf 28} 4415--21

\bibitem{GPS}
Gurevich D I, Pyatov P N  and Saponov P A  
1997 {\it Lett. Math. Phys.} {\bf 41} 255--64 

\bibitem{PS2} Pyatov P N and Saponov P A 1996 
Newton relations for quantum matrix algebras of RTT-type
{\it Preprint} IHEP 96-76

\bibitem{IOPS} 
Isaev A, Ogievetsky O, Pyatov P and Saponov P 1997 
Characteristic polynomials for Quantum Matrices
{\it Preprint} CPT-97/P3471 

\bibitem{IOP}Isaev A, Ogievetsky O and Pyatov P 1998
Generalized Cayley-Hamilton-Newton identities 
{\it Preprint}  math.QA/9809047

\bibitem{FRT} Reshetikhin N, Takhtadjan L and Faddeev L 
1990 {\it Leningrad Math. J.} {\bf 1} 193--225

\bibitem{G} Gurevich D I  
1991 {\it Leningrad Math. J.} {\bf 2}  801--28

\bibitem{D} Drinfeld V G 1986 Quantum groups {\it Proc. ICM (Berkeley)} vol~1
(AMS, Providence, 1987) pp~798--820

\bibitem{R}Reshetikhin N Yu 1990 
{\it Lett. Math. Phys.} {\bf 20} (1990) 331

\bibitem{FM} Freidel L and Maillet J M 1991 {\it Phys. Lett.} B {\bf 262} 278--84

\bibitem{KSk} Kulish P P and Sklyanin E K 1992
{\it J. Phys. A: Math. Gen.} {\bf 25} 5963--75

\bibitem{M} Maillet J M 1990
{\it Phys. Lett.} B {\bf 245}  480-6

\end{thebibliography}
\end{document}